\documentclass{amsart}
\usepackage{a4}
\usepackage{amssymb}
\usepackage{lcy}
\usepackage{verbatim}

\def\wid{\operatorname{wid}}
\def\mba{\mbox{\boldmath{$a$}}}
\def\mbc{\mbox{\boldmath{$c$}}}

\newtheorem{Th}{Theorem}[section] 
\newtheorem{Lem}[Th]{Lemma} \newtheorem{Prop}[Th]{Proposition}
\newtheorem{Claim}[Th]{Claim}

\newtheorem{claim-num}{Claim}

\numberwithin{equation}{section}
\renewcommand{\theequation}{\thesection.\arabic{equation}}

\def\inv{^{-1}}

\def\str#1{\langle#1\rangle}

\def\av#1{\overline{#1}}

\def\N{\mathbf N}
\def\Z{\mathbf Z}

\renewcommand{\le}{\leqslant}
\renewcommand{\ge}{\geqslant}

\def\wid{\operatorname{wid}}
\def\mba{\mbox{\boldmath{$a$}}}
\def\mbc{\mbox{\boldmath{$c$}}}

\begin{document}

\title{The palindromic width of a free product of groups}
\author[Valery Bardakov and Vladimir Tolstykh]{}
\address{Valery Bardakov\\ Institute of Mathematics\\
Siberian Branch Russian Academy of Science\\
630090 Novosibirsk\\
Russia}
\email{bardakov@math.nsc.ru}
\address{Vladimir Tolstykh\\ Department of Mathematics\\ Kemerovo State
University\\ 650043 Kemerovo\\  Russia}
\email{tvlaa@mail.ru}
\maketitle
\begin{center}
\sc Valery Bardakov\footnote{Suppoped
by the RFFI, grant \# 02-01-01118}
and Vladimir Tolstykh\footnote{Supported
by a NATO PC-B grant via The
Scientific and Technical Research Council of Turkey (T\"UBITAK)}
\end{center}

\begin{abstract}
Palindromes are those reduced words of free products
of groups that coincide with their reverse words.
We prove that a free product of groups $G$ has infinite
palindromic width, provided that $G$ is not the free
product of two cyclic groups of order two (Theorem \ref{Pal-Width-of-a-Free-Product}). This means
that there is no a uniform bound $k$ such that every
element of $G$ is a product of at most $k$
palindromes. Earlier the similar fact established
for non-abelian free groups. The proof of Theorem
\ref{Pal-Width-of-a-Free-Product} makes use of the
ideas by Rhemtulla developed for the study of the
widths of verbal subgroups of free products.
\end{abstract}

\section*{Introduction}

Let $G$ be a group and $S$ a generating
set of $G$ with $S\inv=S.$ For each $x \in G$ define
the {\it length} $l_S(x)$ of $x$ relative to $S$ to be
the least $k$ such that $x$ is a product of $k$
elements of $S.$ The supremum of the values $l_S(x)$
is called the {\it width} of $G$ with respect to $S$
and is denoted by $\wid(G,S).$ In particular,
$\wid(G,S)$ is either a natural number, or $\infty.$
In the case when $\wid(G,S)$ is a natural number,
every element of $G$ is a product of at most $\wid(G,S)$ elements of $S.$

Generally speaking, the study of the widths of a given
group provides very useful information for
understanding of the structure of the group, for
various combinatorial, algorithmic and model-theoretic
problems concerning the group etc. The concept of
the width provides also a sort of `measure' for generating
sets: the greater the width a particular generating
set $S$ gives, the less `massive' appeared to be $S$
inside the group it generates. With this point
of view such phenomena as, for instance, the finiteness--under
certain natural conditions--of widths of all verbal subgroups of
linear algebraic groups \cite{Mer:VS}, the infiniteness
of widths of verbal subgroups in various
free constructions \cite{Bar:VS,Dobr:VS,Fai:VS, Rhem},
and the finiteness of the width of the symmetric
group of an infinite set relative to any generating set
\cite{Berg} are of considerable importance.

The paper \cite{BST} examines the
primitive and the palindromic widths of a non-abelian
free group. Traditionally, attributes like `primitive'
and `palindromic' refer to the elements of the
corresponding generating sets. Thus, for instance, the
primitive width of a free group is its width
relatively to the set of all primitive elements.
Palindromic words or {palindromes} can
be defined for all free products of groups (in
particular, for free groups) as follows.
Let
\begin{equation} \tag{$*$}
G = \left.\prod_{i \in I}\right.^* G_i
\end{equation}
be a free product of groups. The {\it
palindromes} associated with the free factorization
$(*)$ are those reduced words of $G$ that are
read the same forward and backward. It is clear
that $G$ is generated by palindromes.  Then the
{\it palindromic width} of $G$ is
the width relative to the palindromes.

It is established in \cite{BST} that the palindromic
width of any non-abelian free group is infinite.
Moreover, the paper \cite{BST} contains a sketch of
the proof of the infiniteness of the palindromic width
of any free product of groups at least two of which
are infinite.

In the present paper we deal with arbitrary free
products of groups. It turns out, that almost all free
products have infinite palindromic width; the only
exception is given by the free product of two cyclic
groups of order two, when the palindromic width is
two.

We would like to point out a certain similarity
between the study of the palindromic widths of free
products we are undertaking and the study of the
widths of verbal subgroups of free constructions.
Indeed, like generators of a given verbal subgroup,
{all palindromic words are {\it structurally alike.}
As the reader shall see this makes possible the
application of the methods developed by Rhemtulla
in \cite{Rhem} specially for the study of the widths of
verbal subgroups of free products.

Let $G$ be a free product of groups.  In Section 1 we
consider the case, when one of the free factors of $G$
contains elements of order at least three. We then
use some of the functions introduced by Rhemulla in
\cite{Rhem} to construct a so-called
quasi-homomorphism, say, $\Delta_1$ whose values at
palindromes are bounded by $2.$ A function $\Delta : H
\to \Z,$ where $H$ is a group, is said to be a {\it
quasi-homomorphism}, if for all $x,y \in H$
$$
\Delta(xy) \le\Delta(x)+\Delta(y)+\text{const}.
$$
At the next step we show that
$\Delta_1$ is not bounded from above. It then
follows that $G$ has infinite palindromic width, since
for each $k \in \N$ the values of $\Delta_1$ at words
which are products of at most $k$ palindromes are
uniformly bounded from above.  The definition
of $\Delta_1(g)$ is rather technical and we just
note that the value $\Delta_1(g),$ where $g \in G$
is a reduced word, somehow reflects the information on
occurrences the subwords of the form $\mba\ldots \mba$
and $\mba\inv \ldots \mba\inv$ in the word $g,$ where
$\mba$ is a fixed element of order at least three
from one of the free factors of $G.$

In Section 2 we study the case, when there are no
elements of order greater than two in the free factors
of $G,$ but there is a free factor that has at least
two non-identity elements, say, $\mba$ and $\mbc.$
Then, like in Section 1, we construct a
quasi-homomorphism $\Delta_2$ that vanishes at
palindromes and unbounded from above.  This time, the
value $\Delta_2(g)$ at a given element $g \in G,$
reflects the information on the occurences
of the subwords of the form
$\mba\ldots\mbc\ldots\mba$ in
$g.$

The authors would like to thank Vladimir
Shpilrain, Oleg Belegradek and Oleg
Bogopolsky for helpful discussion.

\section{Case I: elements of order $\ge 3$ in
one of free factors}

Let
\begin{equation}
G = \left.\prod_{i \in I}\right.^* G_i
\end{equation}
be a free product of groups. Any non-identity element $g$
of $G$ can be written as a reduced word
$$
g= v_1\ldots v_n,
$$
where $v_k$ are elements of free factors $G_i$ and for
every $k=1,\ldots,n-1$ the elements $v_k$ and $v_{k+1}$ lie in
different free factors. The elements $v_k$ are said to be the
{\it syllables} of $g.$ Conversely, if $v_1,\ldots,v_n$
are non-trivial elements from the free factors
$G_i$ and for any $k$ the elements $v_k,v_{k+1}$
are members of different factors, then the
product $v_1\ldots v_n$ is a non-identity
element of $G.$ The number of syllables
of a reduced word $g \in G$ is called
the {\it length} of $g.$

One may rewrite the
syllables of $g$ in the reverse order, obtaining
thereby the non-trivial element
$$
\av g= v_n \ldots v_1
$$
of $G.$ We call an element $g$ a {\it palindrome} associated
with the free factorization (\theequation) if
$$
\av g=g.
$$
It is helpful to observe that each palindrome can
be written in the form
\begin{equation} \label{How-Pals-Look-Like}
g = h v \av h
\end{equation}
where $h$ is a reduced word and $v$ a syllable
of $g$; so that $\av h$ represents the
reflection of $h$ relative to the
central syllable $v.$

Clearly, the elements of the free factors $G_i$ participating
in a decomposition
$$
G =\left.\prod_{i \in I}\right.^* G_i
$$
are palindromes. Thus $G$ is generated by palindromes,
and one can define the {\it palindromic width} of $G$
as the width of $G$ relative to the set of
all palindromes.

As we noted in the Introduction, the palindromic width
of a free product is almost always infinite. We begin
therefore with the description of the {\it only}
exceptional case.

\begin{Claim} \label{Z_2*Z_2}
Let $G=A * B$ be a free product
of cyclic groups of order two.
Then the palindromic width
of $G$ is two.
\end{Claim}

\begin{proof}
Let $A=\str a$ and $B=\str b.$ If
a reduced word $g \in G$ begins
with $a,$ then
$$
g =(ab)^k \text{ or } g=(ab)^k a
$$
for a suitable $k \in \N.$ It is easy to see that any
word of the form $(ab)^ka$ is a palindrome. Any word of the form
$(ab)^k =(ab)^{k-1}a \cdot b$ is a product of two palindromes.
\end{proof}

The bulk of this section is devoted to the proof of
the infiniteness of the palindromic width of a free
product of non-trivial groups $G=A*B$ such that one the groups,
say, $A$ has elements of order at least
three. In the next section we consider the case when
both factors $A,B$ of the free product $A*B$ have no
elements of order $\ge 3,$ but one of the factors has
at least three elements. We shall see later that the
general case can be easily reduced to the case of a
free product of two groups.

In both of the described cases we shall be
actually able to prove that the width of $G=A*B$ relative
to some {\it superset} of the palindromes, consisting of
palindromic-like words, is infinite.
Clearly, this will imply the infiniteness of the
palindromic width of $G.$ With this idea in mind, we
introduce the following defintion.

Assume that $C$ is a subset of $A\cup B.$
Consider the alphabet $X = C \cup \{\omega\},$
where the symbol $\omega$ belongs neither to $A,$
nor to $B.$ For every reduced word
$$
g=v_1\ldots v_n
$$
of $G$ we define the word $\alpha(g)$ over the
alphabet $X$ by replacing the syllables of $g$ that are
not in $C$ by $\omega.$ We then call the word $g$ a
{\it $C$-palindrome} if $\alpha(g)$ is a
palindrome over the alphabet $X$ (that is,
$\alpha(g)$ `reads the same backward and forward' as
a word over the alphabet $X$.)
It is worth mentioning that {\it any palindrome is
a $C$-palindrome for every $C \subseteq A\cup B.$}

Applying the idea of representation
\eqref{How-Pals-Look-Like} of palindromes, we see that a typical
reduced $C$-palindrome looks like \begin{equation}
\label{How-C-Pals-Look-Like}
h v \widetilde h,
\end{equation}
where $v$ is a non-trivial element from
$A\cup B,$ $h$ is a reduced word of $G$
and $\widetilde h$ denotes an {\it arbitrary} word
which is, say, the result of a `poor' reflection
of the word $h,$ that is, if
$$
h =v_1\ldots v_k,
$$
then
$$
\widetilde h = v_1'\ldots v_n'
$$
and $v_i \in C$ implies that $v'_{n-i+1}=v_i.$
Let us stress that $\widetilde{\phantom a}$
denotes a one-placed {\it predicate}, not an operation,
like the symbol $\av{\phantom a}$.

\begin{Prop} \label{There-Is-an-Element-of-Order>2}
Let $G=A*B$ be a free product non-trivial
groups such that $A$ contains elements
of order $\ge 3.$ Then the palindromic
width of $G$ is infinite.
\end{Prop}

\begin{proof}

Let us fix a non-identity element $\mba \in A$
of order greater than two. We have that
$\mba\inv \ne \mba.$ We are going to find a
quasi-homomorphism $\Delta_1 : G \to \Z$ that takes
reasonably small values at
$\{\mba,\mba\inv\}$-palindromes (and consequently at
palindromes.)

Let
$$
g = v_1\ldots v_n
$$
be a reduced word of $G.$ Suppose that
there are at least two occurrences
of the aforesaid fixed element $\mba$
in $g.$ If $v_i=\mba$ and $v_j=\mba$
are consecutive occurrences
of $\mba$ in $g,$ then we call
the subword
$$
v_{i+1}\ldots v_{j-1}
$$
by an {\it $\mba$-segment} of $g$ \cite{Rhem}.
The length of this segment, the number
$j-i-1,$ is {\it odd} as the reader
may easily see. The $\mba\inv$-segments
of the reduced words are defined
in a similar fashion.

For example, a word
$$
\mba b_1 \mba b_2 a_1b_3 \mba\inv b_4 \mba b_5 \mba\inv
$$
where $b_i$ are non-identity
elements of $B$ and $a_1 \in A \setminus \{\mba,\mba\inv\}$ contains
two $\mba$-segments, namely,
$$
b_1 \text{ and }  b_2 a_1b_3\mba\inv b_4
$$
of length $1$ and $5$ respectively
and one $\mba\inv$-segment of length $3,$ namely,
$b_4 \mba b_5.$

For each $k \in \N$ we define the following two
functions on the set of all reduced
words of $G$:
$$
d_k(g) = \text{ the number of $\mba$-segments of $g$ of length $2k+1$}
$$
and
$$
d_k^*(g) = \text{ the number of $\mba\inv$-segments of $g$ of length $2k+1$}.
$$
For every $k \in \N$ set also
$$
t_k(g) =d_k(g)-d_k^*(g).
$$
The functions $d_k,d_k^*$ and $t_k$
were introduced in the paper \cite{Rhem}
by Rhemtulla.

Clearly, for every reduced word $g \in G$
$$
d_k(g)=d_k^*(g\inv)
$$
and consequently
$$
t_k(g)+t_k(g\inv)=0.
$$

\begin{Lem}[\mbox{\cite[Lemma 2.11]{Rhem}}] \label{All-But-Nine}
The formula
$$
t_k(gh)=t_k(g)+t_k(h)\quad \forall g,h \in G
$$
holds for all but at most $9$
natural numbers $k.$
\end{Lem}

Now we are ready to construct
a desired quasi-homomorphism. Let $g$ be
a reduced word from $G$ and let
$$
\Delta_1(g) =\sum_{k=0}^\infty r_k,
$$
where $r_k \ge 0$ is the remainder
of division of the number $t_k(g)$
on $2$ (note that only finitely
many terms of the above series
might be non-zero.)

It follows immediately from Lemma \ref{All-But-Nine}
that
\begin{equation}  \label{Delta1-Is-a-Quasi-Hom}
\Delta_1(gh) \le \Delta_1(h) + \Delta_1(h)+9\quad \forall g,h \in G.
\end{equation}
So that $\Delta_1$ is a quasi-homomorphism.

\begin{Lem}
If $g\in G$ is an $\{\mba,\mba\inv\}$-palindrome,
then $\Delta_1(g) \le 2.$
\end{Lem}

\begin{proof}
According to \eqref{How-C-Pals-Look-Like} a reduced
word $g$ which is an $\{\mba,\mba\inv\}$-palindrome
is
$$
h v \widetilde h
$$
for some $h \in G$ and a non-trivial $v$ in
$A\cup B.$ Clearly, for all $k \in \N$
$$
d_k(h) =d_k(\widetilde h).
$$
We claim therefore that for all but at most
two $k \in \N$
$$
d_k(g) =2d_k(h).
$$
Indeed, the `new' $\mba$-segments that
are not in $h$ or in $\widetilde h$
can occur, if
$$
g =u_1 \mba u_2 v \widetilde u_2 \mba \widetilde u_1,
$$
where $v \ne \mba$ and $\mba$ before $v$ indicates the last from
the left occurrence of $\mba$ in $h,$ or,
if
$$
g =u_1 \mba u_2 \mba \widetilde u_2 \mba \widetilde u_1.
$$
Assuming that $|u_2 v \widetilde u_2|=2m+1,$
we have in the first case that
$$
d_m(g) =2d_m(h)+1.
$$
In the second case, letting $k$ denote
the natural number with $|u_2|=2k+1,$ we get
$$
d_k(g)=2d_k(h)+2
$$
In particular, we see that {\it for at
most one natural $l$ the value $d_l(g)$ is odd.} The
similar argument applied to the $\mba\inv$-segments
of $g$ proves that for at most one $r \in \N$ the value of
$d_r^*(g)$ is odd. Thus at most two of the values
$t_k(g)$ might be odd, which completes the proof of
the Lemma.  \end{proof}

Lemma \theLem\ and \eqref{Delta1-Is-a-Quasi-Hom} imply immediately that
for every $g$ which is a product of
at most $k$ $\{\mba,\mba\inv\}$-palindromes
(palindromes)
$$
\Delta_1(g) \le 11k-9.
$$

Now we prove that $\Delta_1$ is not
bounded from above. It suffices to find a
sequence $(g_n)$ of elements of $G$ such that the sequence
$(\Delta_1(g_n))$ is infinitely increasing. For all
natural $n \ge 1$ set
$$
g_n = v w v w^2\ldots v w^n v,
$$
where
$$
v = b\mba \text{ and } w =b\mba^{-1}
$$
and $b$ is a non-trivial element from $B.$

For every $n \ge 2$ we have
$$
g_n = g_{n-1} (b\mba^{-1})^n b\mba=g'_{n-1} \mba (b\mba^{-1})^n b\mba.
$$
Thus $g_n$
gains only one `new' $\mba$-segment of length
$2n+1$ to be added to the $\mba$-segments of $g_{n-1}.$ This leads to
\begin{align} \label{Counting-D's}
& d_0(g_n)=0, \\
& d_1(g_n)=\ldots = d_n(g_n)=1, \nonumber\\
& d_k(g_n)=0\quad \forall k > n \nonumber
\end{align}
for all $n \ge 1.$

For every $n > 1,$ further, the
word $g_n$ has $\mba\inv$-segments only of length
$1$ and $3.$ Hence
\begin{equation}
d_k^*(g_n) =0 \quad \forall k \ge 2
\end{equation}
for all $n \ge 1.$ Combining \eqref{Counting-D's}
and (\theequation), we arrive at the inequality
$$
\Delta_1(g_n)=r_0+r_1+n-1 \ge n-1
$$
valid for all $n \ge 1$ (recall that $r_k$
is the remainder of the division of $t_k(g)$
on $2$.) So that the sequence $(\Delta_1(g_n))$ is
infinitely increasing, as desired.
The Proposition is proven.
\end{proof}

\section{Case II: two non-identity elements in
one of free factors}

We consider now the  remaining case of free products
of two groups.

\begin{Prop} \label{No-Elements-of-Order>2}
Let $G=A*B$ be a free product of groups, where $|A|
\ge 3$ and $|B|\ge 2$ and both $A,B$ contain no
elements of order $\ge 3.$ Then the palindromic width
of $G$ is infinite.  \end{Prop}

\begin{proof}
It follows from the conditions that all non-identity
elements of $A$ and $B$ are of order two.  This
implies that $A$ and $B$ are abelian, and hence $|A|
\ge 4.$ Thus we may pick up two distinct non-identity
elements $\mba$ and $\mbc$ from $A.$

We shall consider $\{\mba,\mbc\}$-palindromes of $G$
and shall prove that the width of $G$ relative to
those is infinite.  This will imply the infiniteness
of the palindromic width of $G.$ We shall once again
exploit the ideas of Rhemtulla's paper \cite{Rhem}.

Since $\mba\inv$-segments are no
longer useful, we work instead
with $\mba$-segments that contain
a {\it unique} occurence of $\mbc$; such
an $\mba$-segment looks like
\begin{equation}
\underbrace{*\ldots *}_m \mbc \underbrace{*\ldots *}_n
\end{equation}
where there are no occurrences
of neither $\mba$, nor $\mbc$
among syllables $*.$ We
call a segment (\theequation)
an {\it $\mba$-segment of type} $(m,n)$ if $m$
is the number of syllables
before, and $n$ is the number
of syllables after $\mbc.$

The following functions on the set
of reduced words of $G$ were introduced by Rhemtulla in \cite{Rhem}:
\begin{align*}
& d_{m,n}(g) = \text{ the number of $\mba$-segments of type $(m,n)$ in $g$},\\
& t_{m,n}(g)=d_{m,n}(g)-d_{n,m}(g),
\end{align*}
where $m,n$ are arbitrary naturals.

\begin{Lem} [\mbox{\cite[p. 581]{Rhem}}]
{\em (i)} For all pairs
of natural numbers $(m,n)$ we have
that the formulae
\begin{align*}
& d_{m,n}(g) = d_{n,m}(g\inv),\\
& t_{m,n}(g)+t_{m,n}(g\inv)=0,
\end{align*}
hold for each reduced word $g$ in $G;$

{\em (ii)} for all but at most $8$
pairs $(m,n)$ of naturals we have
that the formula
$$
t_{m,n}(gh) = t_{m,n}(g)+t_{m,n}(h)
$$
holds for each pair $g,h$ of
reduced words of $G.$
\end{Lem}

This, as above, advices the following
choice of a quasi-homomorphism
to establish the infiniteness of the
palindromic width:
$$
\Delta_2(g) =\sum_{0 \le m < n} t_{m,n}(g).
$$
The mapping $\Delta_2$ is a quasi-homomorphism, since by
Lemma \theLem
$$
\Delta_2(gh) \le \Delta_2(g)+\Delta_2(h)+8
$$
for all $g,h \in G.$

\begin{Lem}
Assume $g$ is an $\{\mba,\mbc\}$-palindrome.
Then $\Delta_2(g) = 0.$
\end{Lem}

\begin{proof}
By the general description \eqref{How-C-Pals-Look-Like}
$g$ is of the form
$$
h v \widetilde h,
$$
where $v$ is a non-identity element of $A\cup B.$
Calculating the value of $\Delta_2$ at $g$
we are interesting only in those $\mba$-segments
of type $(m,n)$ for which $m \ne n.$ Now

1) if an $\mba$-segment of type $(m,n),$ say,
$$
\mba u_1 \mbc u_2 \mba
$$
occurs in $h,$ then its `reflection' in $\widetilde h$ is
$$
\mba \widetilde u_2 \mbc \widetilde u_1 \mba
$$
and hence it is of type $(n,m)$ and vice versa;

2) the case, when the central element
equals $\mba$ may bring a `new' $\mba$-segment
of type $(m,n),$ but also necessarily a `new'
segment of type $(n,m)$;

3) finally, if the central element is $\mbc,$
the only `new' $\mba$-segment of type $(m,n)$
that may occur in such a configuration is
of type $(k,k)$ for some $k \in \N.$

Summing up, we see that for all pairs $(m,n)$ of naturals
with $m < n$
$$
d_{m,n}(g)=d_{n,m}(g).
$$
Hence $t_{m,n}(g)=0$ for all such pairs,
and $\Delta_2(g)=0.$
\end{proof}

Lemma \theLem\ implies that for any
word $g$ which is a product of
at most $k$ $\{\mba,\mbc\}$-palindromes
(palindromes) we have that
$$
\Delta_2(g) \le 8k-8.
$$
Now we are going to find a sequence $(g_n)$ of elements
of $G$ with
$$
\lim_{n \to \infty} \Delta_2(g_n) =\infty.
$$
A non-trivial element $f$ which is
not equal to $\mba,\mbc$ can be
found in $A.$ Take also some non-trivial
$b \in B.$ For all $n \ge 1$ let
$$
g_n =v w v w^2 v\ldots v w^n v,
$$
where $v=b\mba b\mbc$ and $w =bf.$
The element $g_1$ which is equal to
$$
vwv = b\mba b\mbc bf b\mba b\mbc
$$
contains the unique $\mba$-segment
of type $(m,n),$ where $m < n,$
namely, the segment of type $(1,3).$
Thus
$$
\Delta_2(g_1)=1.
$$
Let $n > 1.$ We have
$$
g_n =g_{n-1} w^n v =g'_{n-1}v w_n v=g'_{n-1} b\mba b\mbc (bf)^n b\mba b \mbc.
$$
Then $g_n$ acquires the (only) `new' $\mba$-segment
of type $(1,2n+1)$ to be added to the $\mba$-segments
of $g_{n-1}$ of types $(1,2k+1),$ where
$k=1,\ldots,n-1.$ This proves that
$$
\Delta_2(g_n)=n
$$
for all $n \ge 1.$ The proof of the Proposition
is now completed.
\end{proof}

\begin{Th} \label{Pal-Width-of-a-Free-Product}
Let
\begin{equation} \label{G-as-a-Free-Product}
G = \left.\prod_{i \in I}\right.^* G_i
\end{equation}
be a free product of non-trivial groups. The
palindromic width of $G$ with respect to the
palindromes associated with the free factorization
{\em (\theequation)} is infinite if and only if at
least one of the free factors has more than two elements
or there are at least three free factors. The
palindromic width of the free product of two cyclic
groups of order two is two.
\end{Th}

\begin{proof}
The case of a free product of
two groups is subject of Claim \ref{Z_2*Z_2},
Proposition \ref{There-Is-an-Element-of-Order>2}
and Proposition \ref{No-Elements-of-Order>2}.
Suppose now that $G$ is a free product
of at least three groups. Let $i_0$
be an index from $I.$ Set
$$
A =\left.\prod_{i \ne i_0}\right.^* G_i \text{ and } B=G_{i_0}.
$$
We have that
\begin{equation}
G =A*B
\end{equation}
and then $G$ is a free product of two groups, one
of which, namely, $A$, is infinite.
Then the width of $G$ relative
to the palindromes associated
with the free decomposition (\theequation) is infinite.
On the other hand, the set of all palindromes
associated with \eqref{G-as-a-Free-Product}
is contained in the set of all palindromes
associated with (\theequation). The
width of $G$ with respect to the
former set must be therefore infinite.
\end{proof}


\begin{thebibliography}{99}


\bibitem{Bar:VS} V.~G.~Bardakov, `On the width of
verbal subgroups of some free constructions',
{Algebra i Logika} (5) {36}
(1997), 494--517.

\bibitem{BST} V.~Bardakov, V.~Shpilrain,
V.~Tolstykh, `On the palindromic and
primitive width of a free group', submitted.

\bibitem{Berg} G.~Bergman, `Generating infinite symmetric groups',
preprint, http://math.berkeley.edu/gbergman/papers/.

\bibitem{Dobr:VS} I.~V.~Dobrynina, `On the width in free
products with amalgamation',
{Mat. Zametki} (3) {68} (2000), 353--359.


\bibitem{Fai:VS} V.~Faiziev, `A problem of
expressibility in some amalgamated products of groups',
{J. Austral. Math. Soc.} {71} (2001),
105--115.


\bibitem{Mer:VS} Yu. I. Merzlyakov, `Algebraic linear
groups as full groups of automorphisms and the closure
of their verbal subgroups', {Algebra i Logika} (1)
{\bf 6} (1967), 83--94.

\bibitem{Rhem} A.~H.~Rhemtulla, `A problem of bounded
expressibility in free products', { Proc. Camb.
Phil. Soc.} (3) {\bf 64} (1969), 573--584.


\end{thebibliography}
\end{document}